\documentclass[12pt]{amsart}
\usepackage{amsfonts,amssymb,hyperref}
\newtheorem{thm}{Theorem}[section]
\newtheorem{cor}[thm]{Corollary}
\newtheorem{lem}[thm]{Lemma}
\newtheorem{prop}[thm]{Proposition}
\theoremstyle{definition}

\theoremstyle{remark}

\numberwithin{equation}{section}
\newcommand{\norm}[1]{\left\Vert#1\right\Vert}
\newcommand{\abs}[1]{\left\vert#1\right\vert}

\newcommand{\Real}{\mathbb R}

\newcommand{\To}{\longrightarrow}

\newcommand\Ind{\rm Ind}
\newcommand\Str{\mathop{\rm Str}}
\newcommand\Conf{\mathop{\rm Conf}}
\newcommand\tr{\mathop{\rm Tr}\nolimits}
\begin{document}

\title{Berezin kernels and Analysis on Makarevich Spaces}   %
\author{Jacques Faraut,
 Michael Pevzner}
\address{  J.Faraut : Institut de Math\'ematiques de Jussieu, UMR 7586,
Universit\'e de Paris 6, Pierre-et-Marie-Curie, Case 247, 4 place
Jussieu, F-75252 Paris Cedex} \address{ M.Pevzner : Laboratoire de
Math\'ematiques, UMR 6056 CNRS, Universit\'e de Reims, Moulin de la
Housse, BP 1039, F-51687, Reims}
 \email{faraut@math.jussieu.fr, pevzner@univ-reims.fr}

\thanks{First author partially supported by the European Commission (IHP Network HARP)}%

\dedicatory{Dedicated to Gerrit van Dijk on the occasion of his 65-th birthday}%
\begin{abstract}
Following ideas of van Dijk and Hille we study the link which
exists between maximal degenerate representations and Berezin
kernels.

\bigskip

We consider the conformal group ${\rm Conf}(V)$ of a simple real
Jordan algebra $V$. The maximal degenerate representations $\pi _s$
($s\in {\mathbb C}$) we shall study are induced by a character of a
maximal parabolic subgroup $\bar P$ of ${\rm Conf}(V)$. These
representations $\pi _s$ can be realized on a space $I_s$ of smooth
functions on $V$. There is an invariant bilinear form ${\mathfrak
B}_s$ on the space $I_s$. The problem we consider is to diagonalize
this bilinear form ${\mathfrak B}_s$, with respect to the action of
a symmetric subgroup $G$ of the conformal group ${\rm Conf}(V)$.
This bilinear form can be written as an integral involving the
Berezin kernel $B_{\nu }$, an invariant kernel on the Riemannian
symmetric space $G/K$, which is a Makarevich symmetric space in the
sense of Bertram. Then we can use results by van Dijk and Pevzner
who computed the spherical Fourier transform of $B_{\nu }$. From
these, one deduces that the Berezin kernel satisfies a remarkable
Bernstein identity :
$$D(\nu )B_{\nu } =b(\nu )B_{\nu +1},$$
where $D(\nu )$ is an invariant differential operator on $G/K$ and
$b(\nu )$ is a polynomial. By using this identity we compute a Hua
type integral which gives the normalizing factor for an
intertwining operator from $I_{-s}$ to $I_s$. Furthermore we
obtain the diagonalization of the invariant bilinear form with
respect to the action of the maximal compact group $U$ of the
conformal group ${\rm Conf}(V)$.

\end{abstract}
\maketitle \thispagestyle{empty}
\section{Berezin kernels on Makarevich Spaces}
In this section we shall introduce the notion of a Berezin kernel
on a symmetric space of a particular type by mean of Jordan
algebraic methods.

\subsection{Jordan algebras and their conformal groups}

A finite dimensional vector space $V$ on $\mathbb R$ or $\mathbb C$
is a \emph{Jordan algebra} if it is endowed with a bilinear map
$(x,y)\to xy$ from $V\times V$ into $V$ satisfying the two following
axioms:
\begin{itemize}
\item{(J1)} $xy=yx,\quad\forall x,y\in V$,
\item{(J2)} $x(x^2y)=x^2(xy),\quad x,y\in V$.
\end{itemize}
Let $L(x)\in {\rm End}(V)$ denote for every $x\in V$ the linear map
defined by $L(x)y=xy$ for every $y\in V$. \vskip 10pt

Let $r$ and $n$ denote respectively the rank and the dimension of the
Jordan algebra $V$.
For a regular element $x$, the minimal polynomial $f_x$ is of degree $r$,
$$f_x(\lambda )=\lambda ^r -a_1(x)\lambda ^{r-1}+\cdots +(-1)^ra_r(x).$$
The coefficient $a_j$ is a homogeneous polynomial of degree $j$, $\Delta (x)=a_r(x)$ is the \emph{Jordan polynomial}, and ${\rm tr}\, (x)=a_1(x)$ is the \emph{Jordan trace} of $x$.
(See \cite{[FK]} p.28 for more details).

Define the so-called \emph{quadratic representation} of $V$ by:
$$
P(x)=2L(x)^2-L(x^2),\qquad x\in V.
$$
One checks that $\Delta(x)^{2}=\det(P(x))^{\frac r{n}}$, where
$P(x)$ is seen as an endomorphism of $V$.\vskip 10pt

A real Jordan algebra $V$ is said to be \emph{semi-simple} if the
bilinear form $\tr L(xy)$ is non degenerate on $V$, such an algebra
is unital. Furthermore it is called \emph{Euclidean} if the above
bilinear form is positive definite.

An involutive automorphism $\alpha$ of $V$ is called
\emph{Euclidean} if the bilinear form $\tr L(\alpha(x)y)$ is
positive definite on $V$. For a semi-simple Jordan algebra such a
Euclidean involution always exists. \vskip 10pt

Assume from now on that $V$ is a simple real Jordan algebra.
According to the general construction of Kantor-K\"ocher-Tits one
associates to $V$ a simple Lie
group which can be understood as a group of conformal
transformations of the corresponding Jordan algebra (in a sense we
shall explain).

Let us remind this classical construction. The \emph{structure
group} ${\rm Str}(V)$ of $V$ can be defined as the subgroup of
$GL(V)$ of elements $g$ such that there exists a real number
$\chi(g)$ for which
\begin{equation}\label{chi}
    \Delta(g.x)=\chi(g)\Delta(x),\quad x\in V.
\end{equation}
The map $g\to\chi(g)$ is a character of ${\rm Str}(V)$ which is a
reductive Lie group. \vskip 10pt
The Jordan algebra $V$ can be identified with the abelian group $N$
of its own translations via the map $y\To n_y$ from $V$ to $N$,
where $n_y(x)=x+y,\:\forall x\in V$.
The \emph{conformal group} $\Conf(V)$ (or the Kantor-K\"ocher-Tits
group) of the Jordan algebra $V$ is the group of rational
transformations of $V$ generated by translations, elements in
$\Str(V)$ and the inversion map $j:x\To -x^{-1}$. It is a simple Lie
group.
A transformation $g\in\Conf(V)$ is conformal in the sense that, at
each point $x$, where $g$ is well defined, its differential $(Dg)_x$
belongs to the structure group $\Str(V)$.

The subgroup of all affine conformal transformations
$P=\Str(V)\ltimes N$ is a maximal parabolic subgroup of $\Conf(V)$.
Let $\sigma$ be the involution of $\Conf(V)$ given by
$$
\sigma(g)=j\circ g\circ j,\qquad g\in \Conf(V),
$$
where $j\in \Conf(V)$ is the inversion map. We define $\overline
N=\sigma(N)$ and $\overline P:=\Str(V)\ltimes\overline N$.

From the geometric point of view the subgroup $\overline P$ can be
characterized in the following way:
$$
\overline P=\{g\in\Conf(V)'\:|\:g(0)=0\},
$$
where $\Conf(V)'$ is the subset of $\Conf(V)$ of all conformal
transformations well defined at $0\in V$. It is open and dense in
$\Conf(V)$. Moreover $\Conf(V)'=N\Str(V)\overline N$. The map
$N\times\Str(V)\times\overline N\To\Conf(V)'$ is a diffeomorphism.
We shall refer to this decomposition as to the \emph{Gelfand-Naimark
decomposition} of the conformal group.
Furthermore, for every transformation $g\in\Conf(V)$ which is well
defined at $x\in V$, the transformation $gn_x$ belongs to
$\Conf(V)'$ and its Gelfand-Naimark decomposition is given by :
\begin{equation}\label{gelf-naimark}
gn_x=n_{g.x}(Dg)_x\bar n',
\end{equation}
where $(Dg)_x\in\Str(V)$ is the differential of the conformal map
$x\to g.x$ at $x$ and $\bar n'\in\overline N$ (see \cite{[P]} Prop.
1.4).

\vskip 10pt

The flag variety ${\mathcal M}=\Conf(V)/\overline
P$, which is compact, is the \emph{conformal compactification} of $V$. In
fact the map $x\To(n_x\circ j)P$ gives rise to an embedding
of $V$ into $\mathcal M$ as an open dense subset, and every
transformation in $\Conf(V)$ extends to $\mathcal{M}$.

\vskip 15pt
 The
Euclidean involution $\alpha$ of $V$ introduced above also defines
an involution $\theta$ of  $\Conf(V)$
by:
\begin{equation}\label{theta}
 \theta(g)=\alpha\circ j\circ g\circ
j\circ\alpha.
\end{equation} It turns out that $\theta$ is a
Cartan involution of $\Conf(V)$ (see \cite{[P]} Prop. 1.1). So the
fix points subgroup of $\theta$ : $ U=\Conf(V)^{\theta}$ is a
maximal compact subgroup of $\Conf(V)$.

Let us remind that a simple real Jordan algebra is either a real
form of a simple complex Jordan algebra or a simple complex Jordan
algebra considered as a real one. We conclude this section with the
classification of simple Jordan algebras, together with their
conformal groups and maximal compact subgroups $U$.
\medskip

{\footnotesize\begin{tabular}{l|c|c|c|c}
  \hline & Complex   & Euclidean & Non Euclidean & Non Euclidean  \\
   &     Non split      &  Split           &Split & Non split $(m=2\ell)$\\
   \hline
  $V$ & $Sym_m(\mathbb C)$ & $Sym_m(\mathbb R)$ & $\times$ & $Sym_{2\ell}(\mathbb R)\cap M_{\ell}(\mathbb H)$ \\
   $\Conf(V)$& $Sp_m(\mathbb C)$ & $Sp_m(\mathbb R)$ & $\times$ & $Sp(\ell,\ell)$ \\
  $U$ &  $Sp(m)$& $U(m)$ & $\times$ & $Sp(\ell)\times Sp(\ell)$ \\
  &&&&\\
 \hline
$V$ & $M_m(\mathbb C)$ & $Herm_m(\mathbb C)$ & $M_m(\mathbb R)$ & $ M_{\ell}(\mathbb H)$ \\
   $\Conf(V)$& $SL_{2m}(\mathbb C)$ & $SU(m,m)$ & $SL_{2m}(\mathbb R)$ & $SL_{2\ell}(\mathbb H)$ \\
  $U$ & $SU(2m)$ & $S(U_m\times U_m)$ & $SO(2m)$ & $SO_{2\ell}(\mathbb H)$ \\&&&&\\
 \hline
 $V$ & $Skew_{2m}(\mathbb C)$ & $Herm_m(\mathbb H)$ & $Skew_{2m}(\mathbb R)$ & $\times$ \\
   $\Conf(V)$& $SO_{4m}(\mathbb C)$ & $SO^{*}(4m)$ & $SO(2m,2m)$& $\times$ \\
  $U$ & $SO(4m)$ & $U(2m)$ & $SO_{2m}\times SO_{2m}$ & $\times$ \\&&&&\\
 \hline
 $V$ & $\mathbb C\times\mathbb C^{n-1}$ & $\mathbb R\times\mathbb R^{n-1}$ & $\mathbb R^p\times\mathbb R^{q}$
 & $\mathbb R^n$ \\
   $\Conf(V)$ & $SO_{n+2}(\mathbb C)$ & $SO_o(2,n)$ & $SO_o(p+1,q+1)$  & $SO_o(1,n+1)$ \\
  $U$ & $SO_{n+2}$ & $SO_2\times SO_n$ & $SO_{p+1}\times SO_{q+1}$ & $SO(n+1)$  \\&&&&\\
 \hline
 $V$ & $Herm(3,\mathbb O)_{\mathbb C}$ & $Herm(3,\mathbb O)$ & $Herm(3,\mathbb O_s)$ & $\times$ \\
   $\Conf(V)$& $E_7(\mathbb C)$ & $E_{7(-25)}$ & $E_{7(7)}$ & $\times$ \\
  $U$ & $E_7$ & $E_{6(6)}\times SO(2)$ & $SU(8)$ & $\times$\\
  \hline
  Type & IV&I&II&III\\
  \hline
\end{tabular}}\medskip

We shall refer to the Jordan algebras given in the first and forth
columns as to the \emph{non split case}, and to those of the second
and third columns as to the \emph{split case}. \medskip

\subsection{Makarevich Riemannian Symmetric Spaces}

A Makarevich symmetric space is a reductive symmetric space which
can be realized as an open symmetric orbit in the conformal
compactification $\mathcal M$ of a simple real Jordan algebra $V$.
We refer the reader to \cite{[B1],[B2]} and literature cited there
for a detailed description. We shall concentrate our interest on
Makarevich symmetric spaces carrying an invariant Riemannian
metric.

Let $\alpha$ be as previously a Euclidean involution of the Jordan
algebra $V$ and let
$$
V_0:=\{x\in V\,\vert\,\alpha(x)=x\},\quad V_1:=\{x\in
V\,\vert\,\alpha(x)=-x\},
$$
be the corresponding eigenspaces of $\alpha$ on $V$. Notice that the
set $V_0$ is a Euclidean Jordan algebra, whose dimension and rank
will be denoted by $n_0$ and $r_0$. This fact implies that the
interior $\Omega_0$ of the set $\{x^2\,\vert\,x\in V_0\}$ of
"positive" elements in $V_0$ is a symmetric cone in $V_0$. Notice
that $r=r_0$ in the split case and $r=2r_0$ in the non split case.
The Jordan algebra $W=V_0+iV_1$ is a Euclidean real form of the
complexification $V^{\mathbb C}=V+iV$. We will denote by $\Omega $
the symmetric cone of $W$. If $V$ is a simple Jordan algebra of type
I, II, or III, then $W$ is simple, while, if $V$ is of type IV,
$W\simeq V_0\times V_0$.\medskip

 According to \cite{[B1]} one introduces two groups:
 $$
 G:=\{g\in\Conf(V)\,\vert\,(-\alpha)\circ
 g\circ(-\alpha)=g\}_o,\, {\rm and}\, K:=\{g\in
 G\,\vert\,g.e=e\},
 $$
 where the subscript $_o$ stands as usually for the connected
 component of the identity transformation and $e$ denotes the
 identity element of the Jordan algebra $V$. It follows that $K$
 is a maximal compact subgroup of $G$. Moreover, the associated
 Riemannian symmetric space is a real tube domain:
 \begin{equation}\label{makarevic}
    {\mathcal X}:=G/K=\Omega_0+V_1.
\end{equation}

The set $ {\mathcal X}$ is a Riemannian \emph{Makarevich symmetric space}. We shall
refer to the previous description as to the tube realization of
$\mathcal X$.

\bigskip

Such a Riemannian Makarevich symmetric space can be obtained as a
real form of a Hermitian symmetric space of tube type. The transform
$$x\mapsto u=(x-e)(x+e)^{-1}.$$
maps the symmetric space $\mathcal X$ onto a bounded domain
$\mathcal D$, which is the unit ball in $V$ with respect to a
so-called spectral norm. Its inverse is the \emph{Cayley transform}
$$
c\,:u\ \mapsto x=c(u)=(e+u)(e-u)^{-1}.
$$

\bigskip
If $V$ is a Euclidean Jordan algebra, then $V_0=\{0\}$, and
$\mathcal X$ is a symmetric cone. If $V$ has a complex structure,
then $\mathcal X$ is a Hermitian symmetric tube.

The following table gives the classification of Riemannian
Makarevich symmetric spaces obtained in this way.
\medskip

{\footnotesize
\begin{tabular}{c|c|c|c|c}
  \hline
   & Complex   & Euclidean & Non Euclidean & Non Euclidean  \\
   &           &            &Split & Non split $(m=2\ell)$\\
   \hline & $C_m$ & $A_{m-1}$& $D_m$& $C_{\ell}$\\
  \hline
  $V$ & $Sym_m(\mathbb C)$ & $Sym_m(\mathbb R)$ & $\times$ & $Sym_{2\ell}(\mathbb R)\cap M_{\ell}(\mathbb H)$ \\
  $G$& $Sp(m,{\Real})$ &
$SL(m,{\Real})\times {\Real}_+$ & $\times $ & $Sp(\ell ,{\mathbb C})$  \\
  $K$ &  $U_m$ & $SO_m$ & $\times $ & $Sp(\ell )$ \\&&&&\\
 \hline
$V$ & $M_m(\mathbb C)$ & $Herm_m(\mathbb C)$ & $M_m(\mathbb R)$ & $ M_{\ell}(\mathbb H)$ \\
   $G$& $SU(m,m)$  &  $SL(m,{\mathbb C})\times {\Real}_+$ & $SO(m,m)$ &
$Sp(\ell ,\ell )$ \\
  $K$& $S(U_m\times U_m)$ & $SU_m$ &
$SO_m\times SO_m$ & $Sp(\ell )\times Sp(\ell )$  \\&&&&\\
 \hline
 $V$ & $Skew_{2m}(\mathbb C)$ & $Herm_m(\mathbb H)$ & $Skew_{2m}(\mathbb R)$ & $\times$ \\
   $G$&  $SO^*(4m)$ & $SL(m,{\mathbb H})\times {\Real}_+$ & $SO(2m,{\mathbb
C})$ & $\times $ \\
  $K$ & $U_{2m}$ & $Sp(m)$ & $SO_{2m}$ & $\times $  \\&&&&\\
 \hline
 $V$ & $\mathbb C\times\mathbb C^{n-1}$ & $\mathbb R\times\mathbb R^{n-1}$ & $\mathbb R^p\times\mathbb R^{q}$
 & $\mathbb R^n$ \\
   $G$ &  $SO_0(2,n)$&$SO_0(1,n-1)\times {\Real}_+$&$SO_0(1,p)\times
SO_0(1,q)$& $SO_0(1,n)$\\
  $K$ &   $SO_2\times SO_n$ & $SO(n-1)$ &
$SO(p)\times SO_q$ & $SO_n$ \\&&&&\\
 \hline
 $V$ & $Herm(3,\mathbb O)_{\mathbb C}$ & $Herm(3,\mathbb O)$ & $Herm(3,\mathbb O_s)$ & $\times$ \\
   $G$&$E_{7(-25)}$ & $E_{6(-26)}\times {\Real}_+$ & $SU^*(8)$ &
$\times $  \\
  $K$ &  $E_6\times SO_2$ & $F_4$ & $Sp(4)$ & $\times $ \\
  \hline
\end{tabular}
}
\bigskip

The second row of this table represents the root system of the pair
$(\mathfrak g,\mathfrak a)$, where $\mathfrak g$ is the Lie algebra
of  $G$, and $\mathfrak a$ is a Cartan subspace.
\bigskip

\noindent
\emph{Example 1}.
Let $V=M(m,{\mathbb R})$, the space of $m\times m$ real matrices, be equipped
with the Jordan product $x\circ y=\frac 12 (xy+yx).$ Then the
Jordan determinant  coincides with the usual matrix determinant :
$\Delta (x)=\det x.$
The structure group ${\rm Str}(V)$ is the group
$S\bigl(GL(m,{\Real})\times GL(m,{\Real})\bigr)$, acting on $V$ by
$$x\mapsto g_1xg_2^{-1} \quad \bigl(g_1,g_2\in GL(m,{\Real})\bigr).$$
The conformal group ${\rm Conf}(V)$ is the group
$SL(2m,{\Real})/\{\pm I\}$ acting on $V$ by
$$x\mapsto (ax+b)(cx+d)^{-1} \quad
{\rm if}\ g=\left(%
\begin{array}{cc}
  a & b \\
  c & d \\
\end{array}%
\right).$$

The differential of a conformal transformation $g$ is given by $
\left(Dg\right)_xy=h_1(x)yh_2(x)$, where, under the condition that
$\det c\neq0$,
$$
h_1(x)=(ac^{-1}d-b)(cx+d)^{-1}c,\:h_2(x)=(cx+d)^{-1}.
$$
Since $\det(ac^{-1}d-b)\det c=\det\left(%
\begin{array}{cc}
  a & b \\
  c & d \\
\end{array}%
\right)=1$ we finally get:
$$
\chi\left(\left(Dg\right)_x\right)=\det(cx+d)^{-2}.
$$
The Euclidean involution $\alpha$ on $V=M(m,\Real)$ is given by
the usual matrix transposition : $ \alpha(x)=x^T$, and
$$V_0=Sym _m(\Real ),\ V_1=Skew _m(\Real ).$$
Remind
that an element $g\in\Conf(V)$ belongs to the group $G$ if and
only if
$$
g\circ(-\alpha)=(-\alpha)\circ g.
$$
If $g.x=(ax+b)(cx+d)^{-1}$, the above condition on $g$ to be in $G$
leads to
$$
(xc^T+d^T)(-ax+b)=-(xa^T+b^T)(-cx+d),
$$
or equivalently,
$$
a^Tc+c^Ta=0,\quad b^Td+d^Tb=0,
$$
and, for every $x\in V$,
$$ x(a^T+dc^Tb)=(d^Ta+b^Tc)x.$$
By Schur's
lemma this last condition says that there is $\lambda\in\Real$ such
that $ a^Td+c^Tb=\lambda{\rm Id}_m. $

Define the matrix
$$
\Upsilon=\left(%
\begin{array}{cc}
  0 & {\rm Id}_m \\
  {\rm Id}_m & 0 \\
\end{array}%
\right).
$$
Then the element $g$ belongs to $G$ if and only if $ g^T\Upsilon
g=\lambda \Upsilon.$ Since $\det g=1$, it follows that
$\lambda^{2m}=1,$ and since $G$ is connected, $\lambda=1$.\medskip

Finally we have shown that $G=SO_o(\Upsilon)/\pm{\rm Id}$, where
$SO(\Upsilon)$ stands for the special orthogonal group of the
quadratic form defined by the matrix $\Upsilon$ which has
signature $(m,m)$. Therefore $G\simeq SO_o(m,m)/\pm{\rm Id}$.

\bigskip
 The conformal compactification $\mathcal M$ of $V$ is the
Grassmann manifold $\mathcal M={\rm Gr}(2m,m)$ of $m$-dimensional
vector subspaces in ${\Real}^{2m}$.
\bigskip

 The corresponding Riemannian
Makarevich symmetric space $\mathcal X$ in its tube realization is
the set of matrices $x\in M(m,{\Real })$ with a positive definite
symmetric part:
$$x+x^T\gg 0.$$

Its bounded realization is given by:
$${\mathcal
D}=\{x\in M(m,{\Real})\mid \|x\|_{\rm op}<1\}.$$

\medskip
\noindent \emph{ Example 2.} Let $V$ be the space ${\Real}\times
{\Real}^{n-1}$. One writes $x=(x_0,x_1),\ x_0\in {\Real},\ x_1\in
{\Real}^{n-1}$. Then the Jordan product is defined as follows :
$z=x\circ y$ if $z_0=x_0y_0-(x_1|y_1),\ z_1=x_0y_1+y_0x_1,$ where
$(x_1|y_1)$ is the usual inner product on $\Real^{n-1}$. The Jordan
determinant is then given by $\Delta (x)=x_0^2+\|x_1\|^2.$ (This
Jordan algebra is actually a Jordan field.)

The corresponding structure group is the group ${\rm
Str}(V)={\Real}_+\times O(n)$, and the conformal group ${\rm
Conf}(V)$ is equal to $O(1,n+1)$. The conformal compactification of
$V$ is the $n$-sphere $S^n$.

Consider the Euclidian involution $\alpha : (x_0,x_1)\mapsto
(x_0,-x_1)$. The corresponding groups $G$ and $K$ are isomorphic
respectively to $SO_0(1,n)$ and $SO(n)$.

The Riemannian symmetric space $\mathcal X$ is the real hyperbolic
space of dimension $n$, realized as a half-space
$${\mathcal X}=\{(x_0,x_1)\in {\Real}\times {\Real}^{n-1}
\mid x_0>0\}.$$ The bounded realization of this symmetric space is
isomorphic to the Euclidean ball:
$${\mathcal
D}=\{x\in {\Real}^n \mid \|x\|<1\}.$$ \medskip

\subsection{Iwasawa decomposition and spherical Fourier Transform}
According to (\ref{theta}) the Cartan involution  of the group $G$
is given by
$$
\theta(g)=(-j)\circ g\circ(-j),\quad g\in G.
$$
The Lie algebra $\mathfrak g$, which is an algebra of quadratic
vector fields on $V$, decomposes under $d\theta$ into a direct sum
of two eigenspaces:
$$
\mathfrak g=\mathfrak k+\mathfrak p,
$$
where the Lie algebra $\mathfrak k$ is the Lie algebra of the group
$K$ introduced above. Let us fix a Jordan frame
$\{c_1,...,c_{r_0}\}$ in $V_0$. Then the space $\mathfrak a$ of the
linear vector fields
$$\xi (x)=\sum _{j=1}^{r_0}t_jL(c_j)x\quad (t_j\in {\Real}),$$
is a Cartan subspace in $\mathfrak p$. The root system of the pair
$({\mathfrak g},{\mathfrak a})$ is of type $A$, $C$, or $D$. We
choose the Weyl chamber ${\mathfrak a}^+$ defined by $t_1<\cdots
<t_{r_0}$ in case of type $A$, or $0<t_1<\cdots <t_{r_0}$ in case
of type $C$ or $D$.

For type $A$ the positive roots are :
$\left\{\frac{t_j-t_i}2,i<j\right\}$, for type $D$ :
$\left\{\frac{t_j-t_i}2,i<j,\:\frac{t_j+t_i}2,i\neq j\right\}$ and
for type $C$: $\left\{\frac{t_j-t_i}2,i<j,\:\frac{t_j+t_i}2,i\neq
j,\:t_i\right\}$.

We define $A=\exp {\mathfrak a}$, and, as usual the nilpotent
subgroup $N$. It can be written $N=N_0\ltimes N_1$, where $N_0$ is
a triangular subgroup in ${\rm Str}(V_0)$, and $N_1$ is the group
of translations
$$x\mapsto x+v\quad {\rm with}\ v\in V_1.$$

The \emph{Iwasawa decomposition} can be described as follows: every
$x\in {\mathcal X}$ can be uniquely written
$$x=na\cdot e=n_0a\cdot e+v,$$
with $a\in A$, $n_0\in N_0$, $v\in V_1$. One writes
$$a=\exp {\mathcal A}(x)\quad {\rm with}\ {\mathcal A}(x)\in {\mathfrak a}.$$
If $\lambda \in {\mathfrak a}^*_{\mathbb C}$, we write, for $\xi
(x)=\sum _jt_j L(c_j)x$,
$$\langle \lambda ,\xi \rangle=\sum _{j=1}^{r_0} \lambda _jt_j.$$
Then
$$e^{\langle \lambda ,{\mathcal A}(x) \rangle}=\Delta _{\lambda }(x_0),$$
where $\Delta _{\lambda }$ is the power function of the symmetric
cone $\Omega _0$ in the Euclidean Jordan algebra $V_0$, and $x_0$
denotes the $V_0$-component of $x$:
$$x=x_0+x_1,\quad x_0\in V_0,\ x_1\in
V_1.$$

\medskip

The \emph{spherical functions} $\varphi _{\lambda }$ for
${\mathcal X}=G/K$ are given by
\begin{equation}\label{sphfunction}
   \varphi _{\lambda }(x)
=\int _K \Delta _{\rho -\lambda }\bigl((k\cdot x)_0\bigr) dk,
\end{equation}
where $dk$ is the normalized Haar measure of $K$, and $\rho $ is
the half-sum of the positive roots:
$$
\langle\rho,\xi\rangle=\frac12{\rm Tr}({\rm ad}\:\xi)_{\vert
Lie(N)}=\frac12{\rm Tr}({\rm ad}\:\xi)_{\vert Lie(N_0)}+\frac12{\rm
Tr}({\rm ad}\:\xi)_{\vert Lie(N_1)}.
$$
One obtains $\rho_j=\frac {d_0}4(2j-r_0-1)+\frac{n_1}{2r_0},$ where
the integer $d_0$ is defined by $n_0=r_0+\frac{d_0}2r_0(r_0-1).$

\medskip

The \emph{spherical Fourier transform} of a $K$-invariant
integrable function $f$ on $\mathcal X$ is defined on $i{\mathfrak
a}^*$ by
$${\mathcal F}f(\lambda )=\int _{\mathcal X} \varphi _{-\lambda }(x)f(x)\mu(dx)
=\int _{\mathcal X}\Delta _{\rho +\lambda }(x_0)f(x)\mu(dx),$$ where
$\mu(dx)=\Delta(x_0)^{-\frac nr}dx$ is a $G$-invariant measure on
$\mathcal X$. Here and elsewhere further the measure $dx$ denotes
the Euclidean measure associated to the Euclidean structure defined
on $V$ by $(x\:\vert\:y)={\rm tr}(x\alpha(y))$.

\subsection{Berezin kernels and the Kantor cross-ratio}

The {\it Kantor cross ratio} of four points of a simple Jordan
algebra $V$ is the rational function
\begin{equation}\label{cross-ratio}
\{x_1,x_2,x_3,x_4\}={\Delta (x_1-x_3)\over \Delta (x_2-x_3)}
:{\Delta (x_1-x_4)\over \Delta (x_2-x_4)}.
\end{equation}
 It is invariant under
conformal transformations and extends to the conformal
compactificaion $\mathcal M$ of $V$.

The invariance by translations is obvious. The invariance under
${\rm Str}(V)$ comes from the fact that $\Delta $ is semi-invariant
under ${\rm Str}(V)$. The invariance under the inverse $j:x\mapsto
-x^{-1}$ follows from the formula (\cite{[FK]}, Lemma X.4.4):
$$
\Delta \bigl(j(x)-j(y)\bigr)={\Delta (y-x)\over \Delta (x)\Delta
(y)}.
$$

Furthermore it is proved that a local transformation which preserves
the cross ratio is the restriction of an element of ${\rm Conf(V)}$
( \emph{cf}. \cite{[Kantor]}, Theorem 6).
\medskip
Being given the Kantor cross-ratio (\ref{cross-ratio}) we shall
introduce the kernel $F$ on $\mathcal X$ by
\begin{eqnarray}
F(x,y)&=&\{x,y,-\alpha (x),-\alpha (y)\}\\
& =&{\Delta \bigl(x+\alpha (x)\bigr)\Delta \bigl(y+\alpha
(y)\bigr) \over \Delta \bigl(x+\alpha (y)\bigr)\Delta
\bigl(y+\alpha (x)\bigr)}.\nonumber
\end{eqnarray}
 This function is invariant under $G$ and
positive. The {\it Berezin kernel} $B_\nu$ can be defined for every
$\nu\in\mathbb C$ as
\begin{equation}\label{berezin}
    B_{\nu }(x,y)=F(x,y)^{\frac{r_0}{r}\nu }=\left\{
\begin{array}{l}
F(x,y)^\nu\quad{\rm in}\:{\rm
the}\:{\rm split}\:{\rm case}\\
F(x,y)^{\frac\nu2}\quad{\rm in}\:{\rm the}\,{\rm non}\,{\rm
split}\,{\rm case} \end{array}\right.
\end{equation}
In virtue of the $G$-invariance of the kernel $B_{\nu }$ the
\emph{Berezin function}
$$
\psi _{\nu }(x) =B_{\nu }(x,e) =\left({\Delta (x_0)\over \Delta
\Bigl({x+e\over 2}\Bigr)^{2}}\right)^{\frac{r_0}{r}\nu }
$$
 is $K$-invariant. In this formula $x_0$ denotes as before the $V_0$-component of
$x$ in the decomposition $V=V_0+V_1,\ x=x_0+x_1.$ \medskip

If $\{c_1,\ldots ,c_{r_0}\}$ is a Jordan frame of $V_0$, and if
$x=\sum _{j=1}^{r_0} e^{t_j}c_j,$ then
$$\psi _{\nu }(x)=\prod _{j=1}^{r_0} \Bigl(\cosh {t_j\over 2}\Bigr)^{-2\nu }.$$

\medskip
\begin{thm}\label{thm1} Assume $\Re\nu>\frac n{r_0}-1$. Then the
function $\psi_{\nu}$ is integrable and its spherical Fourier
transform is given by
\begin{equation}\label{fourierspherik}
 {\mathcal F}\psi _{\nu}(\lambda ) =\frac{P(\lambda,\nu)P(-\lambda,\nu)}{Q(\nu)},
\end{equation}
with
\begin{eqnarray*}
  P(\lambda,\nu) &=& \prod_{j=1}^{r_o}\Gamma\left(\frac 12+\nu-\delta+\lambda_j\right),\quad \delta=\frac{n}{2r_0}\\
  Q(\nu) &=& c\prod_{j=1}^{2r_o}\Gamma\left(\nu-\beta_j\right).
\end{eqnarray*}
where the constant $c$\footnote{We denote by $c$ throughout this
paper different constants depending only on $V$.} and the real
numbes $\beta_j$'s depend on the Jordan algebra $V$. In particular,
for $\nu$ real and $\lambda\in i\mathfrak{a}^*$,
$\bar\lambda=-\lambda$ and
$$
 {\mathcal F}\psi _{\nu}(\lambda ) =\frac{|P(\lambda,\nu)|^2}{Q(\nu)}\geq0,
$$
\end{thm}

Remind that the Gindikin $\Gamma$-function of a symmetric cone
$\Omega$ is defined by
$$
\Gamma_{\Omega}(\nu)=\int_{\Omega}e^{-{\rm
tr}(x)}\Delta(x)^{\nu-\frac nr}dx.
$$
If the cone $\Omega$ is irreducible (the corresponding Euclidean
Jordan algebra $W$ is simple), then
$$
\Gamma_{\Omega}(\nu)=(2\pi)^{\frac{n-r}2}\prod_{j=1}^r\Gamma\left(\nu-(j-1)\frac
d2\right),
$$
where the integer $d$ is defined by $n=r+\frac d2r(r-1)$.

Therefore, for types I and II:
$$
Q(\nu)=c2^{-2r\nu}\Gamma_{\Omega}(2\nu).
$$
For type III:
$$
Q(\nu)=c \Gamma_{\Omega}(\nu).
$$
For type IV:
$$
\Gamma_{\Omega}(\nu)=(\Gamma_{\Omega_0}(\nu))^2,
$$
notice that in this case the Jordan algebra $W$ is not simple
anymore : $W\simeq V_0\times V_0$.

 The above theorem was proved by different methods  in
\cite{[DP],[Z01],[UU]}.

\bigskip

 Let $h(z,w)$ stand for the so-called \emph{canonical polynomial}
defined on $V^{\mathbb C}\times V^{\mathbb C}$ by the following
conditions (see \cite{[FK]} p.262):

-It is holomorphic in the first variable and anti-holomorphic in the
second one.

-For every $g\in\Str(V^{\mathbb C})$ one has $h(gz,w)=h(z,g^*w)$.

-For every $x\in W$ one has $h(x,x)=\Delta(e-x^2)$.\medskip

Recall that the Bergman kernel of the Hermitian symmetric space,
whose $\mathcal{X}$ is a real form, is given by
$h(z,w)^{-\frac{2n}r}$.
\bigskip

\emph{Example 1.} For $V=M(m,\mathbb R)$ the canonical polynomial is
given by $h(x,x)=\det(I-xx^T)$.\medskip

\emph{Example 2.} For $V=\Real\times \Real^{n-1}$ the canonical
polynomial is given by $h(x,x)=(1-\Vert x\Vert^2)^2.$\bigskip

Remind that the Cayley transform
 $$
 c: u\mapsto c(u)=(e+u)(e-u)^{-1},
 $$
 maps the bounded domain $\mathcal{D}$ onto the tube
 $\mathcal{X}=\Omega_0+V_1$.

 Define
$$\tilde F(u,v)=F(c(u),c(v)).$$ Since
$-\alpha(c(u))=c(\alpha(u^{-1}))$ we obtain
$$
\tilde F(u,v):=\{u,v,\alpha(u^{-1}),\alpha(v^{-1})\}={\Delta
\bigl(u-\alpha (u^{-1})\bigr)\Delta \bigl(v-\alpha (v^{-1})\bigr)
\over \Delta \bigl(u-\alpha (v^{-1})\bigr)\Delta \bigl(v-\alpha
(u^{-1})\bigr)}.
$$
Furthermore $\tilde F(u,0)=\lim_{v\to 0}\tilde
F(u,v)=\Delta(\alpha(u^{-1})-u)\Delta(\alpha(u))$, and one can
show that $\tilde F(u,0)=h(u,u)$.
\bigskip

Define similarly
\begin{eqnarray*}
  \tilde B_{\nu}(u,v) &=& B_{\nu}(c(u),c(v)), \\
  \tilde\psi_{\nu}(u) &=& \psi_{\nu}(c(u)).
\end{eqnarray*}
Then
$$
\tilde\psi_{\nu}=\left\{
\begin{array}{l}
h(u,u)^\nu\quad{\rm in}\:{\rm
the}\:{\rm split}\:{\rm case},\\
h(u,u)^{\frac\nu2}\quad{\rm in}\:{\rm the}\:{\rm non}\,{\rm
split}\:{\rm case}.
\end{array}\right.
$$

Consider the integral
$$
I(\nu)=\int_{\mathcal D}\tilde\psi_{\nu}(u)h(u,u)^{-\frac nr}du.
$$
It can be written as
$$
\int_{\mathcal D}h(u,u)^{\frac{r_0}r\nu-\frac nr}du.
$$
Since $h(u,u)^{-\frac nr}du$ is an invariant measure on $\mathcal D$
the last integral equals
$$
 I(\nu)=c\int_{\mathcal X}\psi_{\nu}(x)\mu(dx)=c{\mathcal
F}\psi_{\nu}(\rho).
$$
Therefore we obtain
\begin{cor}
$$
I(\nu)={\rm vol}\:{\mathcal D}\frac{Q(\frac nr)}{P(\rho,\frac
nr)P(-\rho,\frac nr)}\cdot\frac{P(\rho,\nu)P(-\rho,\nu)}{Q(\nu)}.
$$
\end{cor}

Eventually we may notice that according to the Jordan algebra's type
some simplifications in the above formula, due to explicit
expressions of $\Gamma$-factors, are possible, see examples below.
\bigskip

\emph{Example 1.} Let $V=M(m,\Real)$. Then $V_0=Sym_m(\Real),
V_1=Skew_m(\Real)$, the invariant measure on $\mathcal X$ is given
by $\det(x_0)^{-m}dx_0dx_1$, and the Berezin function is
$$
\psi_{\nu}(x)=4^{m\nu}\left[\frac{\det
x_0}{\det(e+x)^2}\right]^{\nu},\:{\rm with}\: \phi_\nu(e)=1.
$$
Its push-forward to $\mathcal D=\{u\in V\:|\:\Vert u\Vert_{\rm
op}<1\}$ equipped with the invariant measure $\det(e-uu^T)^{-m}du$
is given by
$$
\widetilde\psi_\nu(u)=\det(e-uu^T)^\nu.
$$

The spherical Fourier transform of $ \psi_{\nu}$ is the given by
\begin{eqnarray*}
  {\mathcal F}\psi_{\nu}(\lambda) &=& \int_{\omega_0+V_1}\left[\frac{\det
x_0}{\det(e+x)^2}\right]^{\nu}\Delta^o_{\rho-\lambda}(x_0)\Delta^o(x_0)^{-m}dx_0dx_1 \\
   &=& c\frac{\prod_{j=1}^m\Gamma(\nu-\frac{m-1}2+\lambda_j)\Gamma(\nu-\frac{m-1}2-\lambda_j)}
   {4^{-m\nu}\prod_{j=1}^m\Gamma(2\nu-j+1)}.
\end{eqnarray*}
Notice that in this case
\begin{eqnarray*}
Q(\nu)&=&c4^{-m\nu}\Gamma_\Omega(2\nu)=c'\Gamma_{\Omega_0}(\nu)\Gamma_{\Omega_0}(\nu+\frac
12)\\&=&c''\prod_{j=1}^m\Gamma(\nu-\frac{j-1}2)\prod_{j=1}^m\Gamma(\nu+\frac
12-\frac{j-1}2).
\end{eqnarray*}
 Moreover
$$
I(\nu)=\frac{\Gamma_{\Omega_0}(\nu-\frac{m-1}2)}{\Gamma_{\Omega_0}(\nu+\frac
12)}.
$$

\emph{Example 2.} In the case when $V=\Real^n=\Real+\Real^{n-1}$ the
Jordan determinant is given by $\Delta(x)=x_0^2+\Vert x_1\Vert^2$
and the Berezin function is given by
$$
\psi_{\nu}(x)=\left[\frac{4x_0}{(1+x_0)^2+\Vert
x_1\Vert^2}\right]^\nu.
$$
Its spherical Fourier transform is equal to
$$
{\mathcal
F}\psi_{\nu}(\lambda)=c\frac{\Gamma(\nu+\lambda-\frac{n-1}2)\Gamma(\nu-\lambda-\frac{n-1}2)}
{\Gamma(\nu)\Gamma(\nu+\frac 12-\frac{n-1}2)}.
$$
And finally,
$$
I(\nu)=c\frac{\Gamma(\nu-n+1)}{\Gamma\left(\nu+\frac12-\frac{n-1}2\right)}.
$$

\subsection{A Bernstein identity}

\bigskip

An interesting consequence of Theorem (\ref{thm1}) is a Bernstein
identity. Let $D(\nu )$ be the invariant differential operator on
${\mathcal X}=G/K$ whose symbol, {\it i.e.} its image through the
Harish-Chandra isomorphism
$$\gamma : {\mathbb D}({\mathcal X})\rightarrow S({\mathbb C}^{r_0}), \quad
D\mapsto \gamma _D(\lambda ),$$ is
$$\gamma _{D(\nu )}(\lambda )
=\gamma _{\nu }(\lambda ) =\prod _{j=1}^{r_0}(\nu +\frac
12-\delta+\lambda _j)(\nu+\frac 12-\delta-\lambda_j).$$

\begin{cor}\label{bernstein} The Berezin kernel (\ref{berezin}) satisfies the
following identity:
$$
D(\nu)B_{\nu }=b(\nu )B_{\nu +1},
$$
where $b(\nu )$ is the polynomial of degree $2r_0$ given by:
$$
b(\nu)=\frac{Q(\nu+1)}{Q(\nu)}=\prod_{j=1}^{2r_o}(\nu-\beta_j).$$
\end{cor}

In the split case (types I and II):
$$
b(\nu)=\prod_{j=1}^{r_0}\left(\nu-\frac
d4(j-1)\right)\left(\nu+\frac12-\frac d4(j-1)\right).
$$
In case when $V$ is of type III $(r=2r_0)$,
$$
b(\nu)=\prod_{j=1}^{2r_0}\left(\nu-\frac{d}2(j-1)\right),
$$
where the integer $d$ is defined through $n=r+\frac d2r(r-1)$.

In case when $V$ is complex (type IV):
$$
b(\nu)=\prod_{j=1}^{r_0}\left(\nu-\frac{d_0}2(j-1)\right)^2.
$$

 This identity has been established for $V$ complex
by Engli$\check{\rm s}$ \cite{[Englis]}, also in a slightly
different form by Unterberger and Upmeier \cite{[UU]}, and has been
generalized in \cite{[Fa]}.

In general we do not know any explicit expression for the
differential operator $D(\nu )$. However an explicit formula has
been obtained in the special case of a simple complex Jordan
algebra, and for $\nu =\delta $ \cite{[Khlif]}:
$$
   D\left(\delta\right )=\Delta (y)^{1+{n\over
r}} \Delta \Bigl({\partial \over \partial z}\Bigr) \Delta
\Bigl({\partial \over
\partial \bar z}\Bigr) \Delta (y)^{1-{n\over r}},$$ where $z=x+iy$.

Observe that, for $V={\mathbb C}$ we get:
$$D\left(\delta\right)=y^2{\partial ^2\over \partial z\partial \bar
z}$$ is nothing but the Laplace-Beltrami operator of the upper
hyperbolic half-plane.

\bigskip
\emph{Example 1.} When $V=M(m,{\Real})$ the Bernstein polynomial and
the Harish-Chandra symbol are respectively given by
\begin{eqnarray*}
b(\nu)&=&\prod_{j=1}^m\left(\nu-\frac{j-1}2\right)\left(\nu+\frac
12-\frac{j-1}2\right),\\
\gamma_\nu(\lambda)&=&\prod_{j=1}^m\left(\nu-\frac12({m-1})+\lambda_j\right)\left(\nu-\frac12({m-1})-\lambda_j\right).
\end{eqnarray*}

\emph{Example 2.} When $V={\Real}\times {\Real}^{n-1}$ the Bernstein
polynomial and the Harish-Chandra symbol are respectively given by
$$
b(\nu)=\nu\left(\nu-\frac
n2+1\right),\quad\gamma_\nu(\lambda)=\left(\nu+\lambda-\frac{n-1}2\right)\left(\nu-\lambda-\frac{n-1}2\right).$$
\bigskip
\subsection{ Hua Integral}
Let us introduce the compact dual $U/K$ of the symmetric space
$\mathcal X=G/K$ in the complexification $\mathcal X^{\mathbb
C}=G^{\mathbb C}/K^{\mathbb C}$ of the symmetric space $\mathcal X$.

Define $\mathfrak u=\mathfrak k+i\mathfrak p$ and let $U$ be the
analytic subgroup of $G^{\mathbb C}$ with the Lie algebra $\mathfrak
u$. The compact symmetric space $\mathcal Y=U/K$ is isomorphic to
the conformal compactification $\mathcal M$ of the Jordan algebra
$V$.

The set $\mathcal{Y}'$ of invertible elements $y$ in $V^{\mathbb C}$
such that $\bar y=y^{-1}$ is open and dense in $\mathcal{Y}$.

For $\nu=-\kappa,\:(\kappa\in\mathbb N)$ the function
$\psi_{-\kappa}$ extends as a meromorphic function on $\mathcal
X^{\mathbb C}$. For $y=\sum\limits_{j=1}^{r_0}e^{i\theta_j}c_j$ we
have
$$
\psi_{-\kappa}(y)=\prod_{j=1}^{r_0}\left(\cos^2\frac{\theta_j}2\right)^{\kappa}.
$$
This shows that this function is well defined on $\mathcal Y$ for
$\kappa\in \mathbb C$, and bounded for $\Re\kappa\geq0$.

Denote by $\mu_o$ the normalized invariant measure on $\mathcal{Y}$
and define
\begin{equation}\label{hua}
    J(\kappa)=\int_{\mathcal{Y}}\psi_{-\kappa}(y)\mu_o(dy).
\end{equation}

\begin{thm}\label{thmhua}
For $\Re\kappa\geq0$ the following identity holds:
\begin{eqnarray*}
  J(\kappa) &=& c
   \frac{\prod_{j=1}^{2r_o}\Gamma(\kappa+1+\beta_j)}{\prod_{j=1}^{r_o}\Gamma(\kappa+\frac
   12+\delta-\rho_j)\Gamma(\kappa+\frac 12+\delta+\rho_j)}.
\end{eqnarray*}
\end{thm}
\proof Assume first that $\kappa\in\mathbb N$. The Bernstein
identity (\ref{bernstein}) implies that
$$
 D_{-\kappa}\psi_{-\kappa}=b(-\kappa)\psi_{-\kappa+1},
$$

By integrating this identity on $\mathcal{Y}$ we obtain
\begin{eqnarray*}
b(-\kappa)J(\kappa+1)&=&\int_{\mathcal{Y}}D_{-\kappa}\psi_{-\kappa}\mu_o(dy)\\
&=&\int_{\mathcal{Y}}\psi_{-\kappa}(y)D_{-\kappa}1\mu_o(dy),
\end{eqnarray*}
since the operator $D_{-\kappa}$ is self adjoint. The constant term
of the differential operator $D_{-\kappa}$ is given by
$$
D_{-\kappa}1=\gamma_{-\kappa}(\rho).
$$

We finally get
\begin{equation}\label{recurHua}
    \gamma_{-\kappa}(\rho)J(\kappa)=b(-\kappa)J(\kappa-1).
\end{equation}

On the other hand the sequence
$$
\phi(\kappa)=\frac{\prod_{j=1}^{2r_o}\Gamma(\kappa+1+\beta_j)}
{\prod_{j=1}^{r_o}\Gamma(\kappa+\frac12+\delta-\rho_j)\Gamma(\kappa+\frac12+\delta+\rho_j)}
$$
satisfies the same recursion relation. Indeed,
$$
\phi(\kappa)=\frac{\prod_{j=1}^{2r_o}(\kappa+\beta_j)}{\prod_{j=1}^{r_o}(\kappa+\frac12+\delta-\rho_j)
(\kappa+\frac12+\delta+\rho_j)}\phi(\kappa-1).
$$
Therefore $J(\kappa)=c\phi(\kappa)$ and because of the ``initial"
condition $J(0)=1$, we finally get
$$
J(\kappa)=\frac{\phi(\kappa)}{\phi(0)}.
$$
The functions $\phi(\kappa)$ and $J(\kappa)$ are holomorphic for
$\Re\kappa>0$ and grow polynomially :
\begin{eqnarray*}
  |\phi(\kappa)| &\leq& c_1(1+|\kappa|)^{N_1} \\
  |J(\kappa)|&\leq& J(0)=1.
\end{eqnarray*}
Indeed, $0\leq \psi_{-\kappa}(y)\leq 1$. Therefore, the Carlson's
theorem implies that
$$
J(\kappa)=\frac{\phi(\kappa)}{\phi(0)},
$$
for every $\kappa\in\mathbb C$ such that
$\Re\kappa>0.\qquad\qquad\qquad\qquad\qquad\qquad\qquad\Box$
\bigskip

By writing $y'=c(iv)$, $J(\kappa)$ becomes an integral over $iV$.
One shows that
$$
\tilde\psi_{-\kappa}(iv):=\psi_{-\kappa}(c(iv))=h(v,-v)^{-\frac{r_0}r\kappa}.
$$
Therefore, $$J(\kappa)=\int_Vh(v,-v)^{-\frac{r_0}r\kappa-\frac
nr}.$$ In fact, in this realization an invariant measure $\mu_o$ on
$\mathcal Y$ is given by $h(v,-v)^{-\frac nr}dv$.

 The above integral has been considered and computed by Hua
in \cite{[Hua]}, Chapter II in several special cases.
\section{Berezin kernels and representations}
Following ideas developed in \cite{[vdh],[vdm],[Hille]} we shall now
give an alternative approach to the theory of Berezin kernels
considering intertwining operators for the maximal degenerate series
representations of a conformal group.

\subsection{Maximal degenerate series representations}

The character $\chi$ of the structure group $\Str(V)$ (\emph{cf}.
(\ref{chi})) can be trivially extended to the whole parabolic
subgroup $\overline P$ by $\chi(h\bar n):=\chi(h)$ for every
$h\in\Str(V),\:\bar n\in\overline N$.

For every $s\in\mathbb C$ we define a character $\chi_{s}$ of
$\overline P$ by
$$
\chi_{s}(\bar p):=|\chi_(\bar p)|^{s+\frac{n}{2r}}.
$$

The induced representation $\widetilde{\pi_{s}}=\Ind_{\overline
P}^{\Conf(V)}\left(\chi_{s}\right)$ of the conformal group acts on
the space
$$
\widetilde I_{s}:=\{f\in C^{\infty}(\Conf(V))\:|\:f(h\bar
p)=\chi_{s}(\bar p)f(h),\forall h\in\Str(V),\bar p\in\overline P\},
$$
by left translations. A pre-Hilbert structure on $\widetilde I_{s}$
is given by
$$\norm f^2=\int_U\abs{ f(u)}^2du,$$ where $U$ is
the maximal compact subgroup of $\Conf(V)$ associated with the
Cartan involution $\theta$, and $du$ is the normalized Haar measure
of $U$.

According to the Gelfand-Naimark decomposition a function
$f\in\widetilde I_{s}$ is determined by its restriction
$f_V(x)=f(n_x)$ on $N\simeq V$. Let $ I_{s}$ be the subspace of
$C^{\infty}(V)$ of functions $f_V$ with $f\in \widetilde I_{s}$. The
conformal group acts on $I_{s}$ by:
\begin{equation}\label{action}
    \pi_{s}(g)f(x)=|A(g,x)|^{s+{n\over 2r}}f(g^{-1}.x),\:g\in\Conf(V),\:x\in V,
\end{equation}
where
$$A(g,x):=\chi_{s}\bigl((Dg^{-1})_x\bigr).$$
This action is usually called the \emph{even maximal degenerated
series} representation of $\Conf(V)$.

One shows that the norm of a function $f(n_x)=f_V(x)\in I_{s}$ is
given by:
\begin{equation}\label{norm}
    \norm f^2=\int_V\vert f_V(x)\vert^2h(x,-x)^{2\Re s}dx,
\end{equation}

The formula (\ref{norm}) implies that for $\Re s=0$ the space
$I_{s}$ is contained in $L^2(V)$ and the representation $\pi_{s}$
extends as a unitary representation on $L^2(V)$.

\bigskip

In order to address the question of irreducibility and unitarity of
these representations we refer the reader to \cite {[P]},
\cite{[S2],[S1]} and to \cite{[Z02]}. \vskip 15pt

 Following the standard procedure
we introduce an intertwining operator between $\pi_{s}$ and
$\pi_{-s}$.

Consider the map $\widetilde A _{s}$ defined on $\widetilde I_{s}$
by
\begin{equation}\label{intertwiner}
    f\To(\widetilde A_sf)(g):=\int_Nf(gjn)dn,\quad\forall g\in\Conf(V),
\end{equation}
where $dn$ is a left invariant Haar measure on $N$. We will see that
this integral converges for $\Re s>{n\over 2r_0}$.


\begin{prop}
For every $f\in\widetilde I_{s}$ the function $\widetilde A_sf$
belongs to $\widetilde I_{-s}$ and the map $\widetilde A_s$ given by
(\ref{intertwiner}) intertwines the corresponding representations of
the conformal group:
\begin{equation}\label{intertw-repr}
    \widetilde{\pi}_{-s}(g)(\widetilde A_sf)=\widetilde A_s(\widetilde{\pi}_{s}(g)f),\:\forall
    f\in\widetilde I_{s},\:g\in\Conf(V).
\end{equation}
\end{prop}

\emph{Proof.} The map $\widetilde A_s$ obviously commutes with the
left action of $\Conf(V)$. We have to check that $\widetilde A_sf$
transforms in an appropriate way under the right action of
$\overline P$. Notice first that $\widetilde A_sf$ is $\overline N$
right invariant. Indeed
\begin{eqnarray*}
(\widetilde A_sf)(g\bar n)&=&\int_Nf(g\bar
njn)dn=\int_Nf(gjn'jjn)dn\\&=&\int_Vf(gjn'')dn''=(\widetilde
A_sf)(g).
\end{eqnarray*}
For what concerns the action of $\Str(V)$ we have for every
$f\in\widetilde I_{s}$:
\begin{equation}\label{AfH}
    (\widetilde A_sf)(g\ell)=\int_Nf(g\ell jn)dn=\int_Nf(gj(\ell^t)^{-1}n)dn,
\end{equation}
where $h^t$ denotes the transformation adjoint to $h\in\Str(V)$ with
respect to the bilinear form $\tr L(xy)$ on $V$.

Indeed, one shows, (see \cite{[FK]}Prop.VIII.2.5) that $j\ell
j=(\ell^t)^{-1},\:\forall \ell\in\Str(V)$. Thus the equality
(\ref{AfH}) implies that
\begin{eqnarray*}
  (\widetilde A_sf)(g\ell) &=& \int_Nf(gj(\ell^t)^{-1}n\ell^t(\ell^t)^{-1})dn \\
  &=& \abs{\det \ell}^{\frac nr}\int_Nf(gjn'(\ell^t)^{-1})dn' \\
   &=& |\det \ell|^{\frac nr-s\frac{n}{2r}}(\widetilde A_sf)(g)=\chi_{-s}(\ell)(\widetilde A_sf)(g)\qquad\qquad\qquad\Box
\end{eqnarray*}

\vskip 10pt

Define the map $A_s: I_s\rightarrow I_{-s}$ by
$$
A_s(f_V)=(\widetilde A_s f)_V\quad (f\in \widetilde I_s).
$$

\begin{prop}
For every $f\in I_s$ we have,
$$(A_sf)(x)=\int _V |\Delta (x-u)|^{2s-{n\over r}}f(u)du.$$
\end{prop}

\proof By definition of $\widetilde A_s$,
$$
(\widetilde A_sf)(n_x)=\int _V f(n_xjn_v)dv\quad (f\in \widetilde I_s).
$$
According to the Gelfand-Naimark decomposition (\ref{gelf-naimark})
for $g=n_xj$ we have
$$
  n_xjn_v=n_{(n_x\circ j).v}(D(n_x\circ j))_v\bar
  n')=n_{(x-v^{-1})}P(v^{-1})\bar n'.
$$

Thus
$$
\left(A_sf\right)(x) =\int_Vf(x-v^{-1})\chi_{s}(P(v)^{-1})dv.
$$
Let $u=x-v^{-1}$. The Jacobian of this transformation being equal to
$|\Delta (v)|^{\frac{2n}r},$ we finally get:
\begin{equation}\label{intertwiner_V}
 \left(A_sf\right)(x)=\int_Vf(u)|\Delta(x-u)|^{2s-\frac{n}r}du.\qquad\quad\qquad\Box
\end{equation}
\medskip
The bilinear form on $ \widetilde I_{-s}\times \widetilde I_{s}$
defined by
$$
\langle f_1,f_2\rangle=\int_Uf_1(u)f_2(u)du
$$
is $\Conf(V)$-invariant. By using the fact that
$$
\langle f_1,f_2\rangle =c\int _V (f_1)_V(x)(f_2)_V(x)dx.
$$
it gives rise to a $\Conf(V)$-invariant bilinear form on $I_s\times
I_s$:

\begin{prop}\label{Bs}
For $\Re s>{n\over 2r_0}$, the bilinear form $\mathfrak{B}_s$,
$$
 \mathfrak{B}_s(f_1,f_2)
=\int \int _{V\times V}|\Delta (x-y)|^{2s-{n\over
r}}f_1(x)f_2(y)dxdy,
 $$
is well defined on $I_s\times I_s$ and is ${\rm Conf}(V)$-invariant.
\end{prop}
\proof The invariance follows from the intertwining property of
$A_s$. It can also be shown by a straightforward verification.
Indeed, let us recall following property of conformal
transformations (\emph{cf}. \cite{[FG]}, Lemma 6.6) :
\begin{equation}\label{FG}
    \Delta(g.x-g.y)^2=A(g,x)A(g,y)\Delta(x-y)^2,\quad g\in\Conf(V),\:x,y\in V.
\end{equation}
Therefore,
\begin{eqnarray*}
  &&\mathfrak{B}_s(\pi_s(g)f_1,\pi_s(g)f_2)= \\
  &=& \iint\limits_{V\times V}|\Delta(x-y)|^{2s-\frac nr}|A(g^{-1},x)|^{s+\frac n{2r}}
| A(g^{-1},y)|^{s+\frac n{2r}} f_1(g^{-1}.x)f_2(g^{-1}.y)dxdy
\end{eqnarray*}
Denote $x=g.x',\:y=g.y'$. Thus $dx=A(g,x')^{\frac nr}dx'$ and
$dy=A(g,y')^{\frac nr}dy'$. Notice that $g\cdot g^{-1}={\rm Id}$ and
therefore $\left(Dg^{-1}\right)_{g.x}\circ\left(Dg\right)_x={\rm
Id}$. It implies that $A(g^{-1},g.x)A(g,x)=1$.

Taking into account the identity (\ref{FG}) we finally get
\begin{eqnarray*}
 &&\mathfrak{B}_s(\pi_s(g)f_1,\pi_s(g)f_2) = \\
  &=& \iint\limits_{V\times V}|\Delta(x'-y')|^{2s-\frac nr}|A(g^{-1},x')|^{s-\frac n{2r}}
  |A(g^{-1},y')|^{s-\frac n{2r}}\\
  &\cdot&A(g,x')^{-s-\frac n{2r}} A(g,y')^{-s-\frac n{2r}}
  |A(g,x')|^{\frac nr}|A(g,y')|^{\frac nr}f_1(x')f_2(y')dx'dy'\\&=&\mathfrak{B}_s(f_1,f_2).
\qquad\qquad\qquad\qquad\qquad\qquad\qquad\qquad\qquad\qquad\Box
\end{eqnarray*}

\bigskip

 According to the general theory (see for instance \cite{[K]})
these maximal degenerate series representations $\pi_{s}$ are
spherical. We shall determine the corresponding $U$-fixed vector of
the representation $\pi_{s}$.

Using the Iwasawa decomposition of $\Conf(V)$ we write
$g=u\ell(g)\bar n$, where $\ell(g)\in\Str(V)$, which is defined up
to a multiplication by an element of $U\cap\Str(V)$ on the left.
Notice by the way that with the above notation we have
$\theta(g^{-1})g=n^{-1}\ell(g)^2\bar n$.

Consider $n_x\in N$ admitting the decomposition $n_x=u\ell(n_x)\bar
n$. Therefore
$$
\ell(n_x)^2=D(\theta(n_{-x})n_x)(0)=P(x^{-1}+\alpha(x))^{-1}P(x)^{-1}.
$$
Let $f^{o}\in\widetilde I_{s}$ be a $U$-invariant function, with
$f^{o}(e)=1$. Then
$$
f^{o}(g)=f^{o}(u\ell(g)\bar
n)=\chi_{s}(\ell(g))f^{o}(u)=\chi_{s}(\ell(g)).
$$
Eventually we get that the $U$-fixed vector $f^{o}_V$ in $I_{s}$ is
given by
$$
f^{o}_V(x)=|\Delta(x^{-1}+\alpha(x))\Delta(x)|^{-(s+\frac{n}{2r})}=h(x,-x)^{-(s+\frac
n{2r})}.
$$

\bigskip

In order to conclude this section we shall discuss conditions under
which the intertwining operator $A_s$ is well defined.

Let us evaluate $A_sf^o$ at the identity of the conformal group.
$$
\widetilde
A_sf^o(e)=\int_Nf^o(jn)dn=\int_Nf^o(n)dn=\int_Vh(x,-x)^{-(s+\frac
n{2r})}dx.$$ Notice that the value $C(s)$ of this integral equals
$$C(s)=J\left(\frac r{r_0}(s-\frac{n}{2r})\right),$$
where $J(\nu)$ is given by (\ref{hua}). For $u\in U$,
$f^o(u)\not=0$. therefore, for $f\in \widetilde I_s$,
$$|f(g)|\leq M|f^o(g)|,$$
with
$$M= \sup _{u\in U}{|f(u)|\over {|f^0(u)|}}.$$
Now it is clear that the integral defining $A_s$ is well defined for
$\Re s>{n\over 2r_0}$. One shows that the map $s\to A_s$ can be
actually extended to whole complex plane as a meromorphic function.
Because of the intertwining property we eventually get:
$$
A_s\circ A_{-s}=C(s)C(-s){\rm Id}_{ I_s}.
$$

As it has been noticed by several authors \cite{[BLIN],[KS],[Z02]}
for some particular values of the parameter $s$ the intertwining
operator $A_s$ turns out to be a differential operator, which,
furthermore, satisfies a generalized Capelli identity. Indeed, the
map $A_s$ may be seen as a convolution operator with the kernel
$|\Delta(x-y)|^{2s-\frac nr}$ which has a meromorphic continuation
with respect to $s$. One shows that, for $2s-\frac
nr=-2\ell,\:\ell\in\mathbb N$, this is a differential operator
proportional to $\Delta\left(\frac{\partial}{\partial
x}\right)^{2\ell}.$
\bigskip

\subsection{Restriction of $\pi_s$ and Canonical representations}
\bigskip

We shall study the representations $\pi_{s}$  when restricted to the
subgroup $G$. The representation space $\widetilde I_s$ can be seen
as a line bundle over the conformal compactification $\mathcal{M}.$

The Makarevich space $\mathcal{X}$ is one of the open orbits of $G$
acting on $\mathcal{M}$. The space $I_{s}({\mathcal X})$ of
functions in $I(s)=I_{s}(V)$ supported in the closure
$\overline{\mathcal X}\subset V$ is then invariant. We shall
consider the corresponding representation $T_s$ of $G$ on
$I_{s}({\mathcal X})$.

The problem is to determine for which values of the parameter $s$
the representation $T_s$ of the group $G$ is unitarizable, and then
to decompose it into irreducible ones. According to an established
terminology taking its roots in \cite{[GVG]} one calls such
representations \emph{canonical representations} of the group $G$.

A key observation is the connection between the canonical
representation $T_s$ and the Berezin kernel. We shall make this link
clear.\medskip

Consider the bilinear form on $I_{s}({\mathcal X})\times
I_{s}({\mathcal X})$ given by
\begin{equation}\label{bformmakar}
    \widetilde{\mathfrak B}_s^{\alpha}(f_1,f_2):=\mathfrak{B}_s(f_1,f_2\circ(-\alpha)),
\end{equation}
where $\mathfrak{B}_s(f_1,f_2)$ is the bilinear form on
$I_{s}(V)\times I_{s}(V)$ introduced in Proposition (\ref{Bs}).
\begin{prop}
The bilinear form $\widetilde{\mathfrak B}_s^{\alpha}$ is invariant
under the action of the group $G$.
\end{prop}
The proof uses the invariance property of the cross-ratio, and the
two following lemm\ae.
\begin{lem}\label{lemAalpha}
For every $g\in G$ the following identity holds:
$$
A(g,-\alpha(x))=A(g,x).
$$
\end{lem}
\proof By definition the transformation $g$ commutes with $-\alpha$
for every $g\in G$: $ g\circ(-\alpha)=(-\alpha)\circ g. $
Differentiating this equality we have
$$
\left(Dg\right)_{-\alpha(x)}\circ(-\alpha)=(-\alpha)\circ
\left(Dg\right)_x,
$$
or equivalently $\left(Dg\right)_{-\alpha(x)}=(-\alpha)\circ
\left(Dg\right)_x\circ(-\alpha).$
 Notice that for every $h\in\Str(V)$  one has
 $$
 (-\alpha)\circ h\circ(-\alpha)=\alpha\circ h\circ\alpha\in\Str(V),
 $$
 and on the other hand $\chi(\alpha\circ h\circ\alpha)=\chi(h)$. Since
$\Delta(\alpha(x))=\Delta(x)$ we finally get
$$
\Delta(\alpha\circ h\circ\alpha(x))=\chi(h)\Delta(x),
$$
which implies the statement for
$h=\left(Dg\right)_x.\qquad\qquad\qquad\Box$
\begin{lem}\label{lemgx0}
Let $x$ be an element in $V$, then for every $g\in G$ the following
identity holds
$$
\Delta((gx)_0)=A(g,x)\Delta(x_0).
$$
(Recall that the subscript $_0$ means the $V_0$-component.)
\end{lem}
Notice that
$\Delta((gx)_0)=\Delta(g.x-g.(-\alpha(x))=\Delta(g.x+\alpha(g.x))$
for every $g\in G$. Let us substitute $y$ by $-\alpha(x)$ in
(\ref{FG}), then
$$
\Delta(g.x-g.(-\alpha(x)))^2=A(g,x)A(g,-\alpha(x))\Delta(x+\alpha(x))^2.
$$
Therefore, according to the lemma (\ref{lemAalpha}), we have
$$
\Delta((gx)_0)=\pm A(g,x)\Delta(x_0).
$$
Since for $g=Id$ the sign in the above equation is positive, and the
group $G$ being connected it remains positive for every $g\in
G.\quad\quad\Box$
\bigskip

Let us introduce the multiplication operator on the space
$I_s(\mathcal X)$ by
$$
M_{s }f(x)=\Delta (x_0)^{s+\frac n{2r} } f(x).
$$
These operators intertwine the canonical representation $T_s$ and
the left regular action $L$ of the group $G$,
$$
M_s\circ T_s(g)=L(g)\circ M_s,\quad g\in G,
$$
where $ \bigl(L(g)f\bigr)(x)=f(g^{-1}\cdot x).$

\begin{prop}
Let us define the bilinear form ${\mathfrak B}^{\alpha}_s$ on
$\mathcal C_c(\mathcal X)\times \mathcal C_c(\mathcal X)$ by
$$
\mathfrak B_s^{\alpha }(F_1,F_1) =\mathfrak B_s^{\alpha
}(M_sf_1,M_sf_2) =\widetilde{ \mathfrak B}_s^{\alpha }(f_1,f_2).
$$
Then
\begin{equation}\label{berezinmakarintegral}
    {\mathfrak B}^{\alpha}_s(F_1,F_2)
=\iint\limits_{\mathcal X\times\mathcal X}{B}_{\nu
}(x,y)F_1(x)F_2(y)\mu(dx)\mu(dy),
\end{equation}
where $B_{\nu}(x,y)$ is the Berezin kernel introduced in
(\ref{berezin}) and \newline $\nu=-\frac{r}{r_0}\left(s-\frac
n{2r}\right)$.
\end{prop}

Therefore the decomposition of the canonical representations reduces
now to a classical problem in spherical harmonic analysis.

 Let $\Lambda $ be the set of parameters $\lambda \in \mathfrak a _{\mathbb C}^*$
for which the spherical function $\varphi _{\lambda }$  given by
(\ref{sphfunction}) is positive definite. Let $\psi $ be a
continuous function of positive type on $G$ which is
$K$-biinvariant. By the Bochner-Godement theorem there is a unique
bounded positive measure $m $ on $\Lambda $ such that
$$\psi (x)=\int _{\Lambda } \varphi _{\lambda }(x)m (d\lambda ).$$
If $\psi $ is integrable, then the measure $m $ is absolutely
continuous with respect to the Plancherel measure, with a density
given by the Fourier transform of $\psi $:
$$m(d\lambda )=\mathcal F \psi (\lambda ){d\lambda \over |c(\lambda )|},$$
where $c(\lambda)$ is the Harish-Chandra $c$-function of the
symmetric space $\mathcal{X}$.

\bigskip

Therefore the problem is :

(1) Determine the set ${\mathcal W}\subset {\Real}$ of values $\nu $
for which the Berezin kernel ${\mathfrak B}^\alpha_\nu$ (or the
corresponding $K$-biinvariant function $\psi _\nu$) is of positive
type.

(2) For $\nu \in {\mathcal W}$, determine the positive measure $m
_\nu$ on $\Lambda $ such that
$$\psi _{\nu }(x)=\int _{\Lambda } \varphi _{\lambda }(x)m_{\nu}(d\lambda ).$$

These problems have been solved for hyperbolic spaces by van Dijk,
Hille, Pasquale \cite{[vdh],[vdpas]}. Notice however that only real
hyperbolic spaces are Makarevich spaces as introduced above.

The case $G=U(p,q)$ was studied by Hille, and  Neretin
\cite{[Hille],[Ner]}, and the case when $G=SO(p,q)$  by Neretin in
\cite{[Ner]}.

\bigskip

For $\nu>\frac n{r_0}-1$ the function $\psi _\nu$ is integrable and,
for $\lambda\in i\mathfrak a^*$:

$$\mathcal{F}\psi _\nu(\lambda)=\frac{|P(\lambda,\nu)|^2}{Q(\nu)}\geq 0.
$$

It follows that $\psi _\nu$ is of positive type, and
$$
\psi _\nu(x)=\int_{i\mathfrak
a^*}\phi_{\lambda}(x)\frac{|P(\lambda,\nu)|^2}{Q(\nu)}\frac{d\lambda}{|c(\lambda)|^2}.
$$

According to Helgason we define the Fourier transform of an
integrable function $f$ on $\mathcal{X}$. The Iwasawa decomposition
of the group $G$ can be written as $G=NAK$. If $k^{-1}g\in
N\exp(H)K$ with $H\in\mathfrak a$, one writes
$$
H=\mathcal A(x,b),\quad x\in gK,\:b\in kM\in B=K/M,
$$
where $M$ is the centralizer of $A$ in $K$.

The Fourier transform of $f$ is the function defined on $i\mathfrak
a^*\times B$ by
$$
\hat f(\lambda,b)=\int_{\mathcal
X}f(x)e^{\langle-\lambda+\rho,\mathcal{A}(x,b)}\rangle dx.
$$

We consider $\hat f(\lambda):=\hat f(\lambda,\cdot)$ as an element
of the space $\mathcal{H}_{\lambda}\simeq L^2(B)$ which carries a
unitary spherical principal series $\Pi _\lambda$ of $G$. Then the
map $f\mapsto\hat f(\lambda)$ intertwines the left regular
representation $L$ and $\Pi _\lambda$.
\begin{thm} 

For $\nu>\frac n{r_0}-1$, the Berezin form is positive definite.
Therefore the representation $T_s$ is unitarizable, decomposes
multiplicity free as a direct integral of spherical principal series
$\Pi _\lambda$ of $G$, which corresponds to the following Plancherel
formula :
$$
\mathfrak{B}_s^{\alpha}(f,\bar f)=\int_{i\mathfrak a^*}\Vert\hat
f(\lambda)\Vert^2\frac{|P(\lambda,\nu)|^2}{Q(\nu)}\frac{d\lambda}{|c(\lambda)|^2},
$$
where $\nu=-\frac r{r_0}(s-\frac n{2r})$ and $f\in C_c(\mathcal X)$.
\end{thm}
\bigskip

\subsection{Berezin kernel on the Riemannian Compact Dual}
In order to study the "deformation" of the regular representation of
the compact group $U$ we shall investigate the spherical Fourier
transform of the Berezin kernels on the compact dual symmetric space
$\mathcal Y$.

 The spherical functions of the symmetric space $\mathcal Y=U/K$ are
given by
\begin{equation}\label{spherY}
    \Phi_{{\bf m}}(x)=\int_K\Delta_{{\bf m}}((k.x)_0)dk=\varphi_{{\bf m}-\rho}(x),
\end{equation}
where the weights ${\bf m}=(m_1,\dots,m_n)\in\mathbb Z^n$ are given
according to different cases by following conditions :
\begin{itemize}
\item If $V$ is a Euclidean Jordan algebra (root system of type A),
then $m_1\geq m_2\geq\dots\geq m_n$.
\item If $V$ is a non Euclidean algebra of split type (root system
of type D), them $-|m|_1\geq m_2\geq\dots\geq m_n$.
\item If $V$ is either a complex or  a non Euclidean Jordan algebra of non split
type (root system of type C), then $0\geq m_1\geq m_2\geq\dots\geq
m_n$.
\end{itemize}
\medskip

The spherical Fourier coefficients of a $K$-invariant integrable
function $f$ on $\mathcal Y$ are given by
$$a({\bf m})=\int _{\mathcal Y} f(y)\Phi _{\bf m}(y)\mu _0(dy).$$
For $\Re \kappa \geq 0$, the Berezin function $\psi _{-\kappa}$ is
bounded on $\mathcal Y$. We will determine its spherical Fourier
coefficients
$$
a_{\kappa}({{\bf m}})=\int_{\mathcal Y}\psi_{-\kappa}(y)\Phi_{{\bf
m}}(y)d\mu_0(dy).
$$
Notice that $a_{\kappa}(0)=J(\kappa)$.
\begin{thm}\label{thmcoefffourier}
The Fourier coefficients of the Berezin kernel function
$\psi_{-\kappa}$ are given by
$$
a_{\kappa}({{\bf
m}})=J(\kappa)\frac{\prod_{j=1}^{n}\Gamma(\kappa+\frac12+\delta-\rho_j)\Gamma(\kappa+\frac12+\delta+\rho_j)}
{\prod_{j=1}^{n}\Gamma(\kappa+\frac12+\delta-\rho_j+m_j)\Gamma(\kappa+\frac12+\delta+\rho_j-m_j)}.
$$
\end{thm}

\emph{Proof} We use the same method as in \cite{[Z03]}.

(a) We show first that the measure on $\mathcal Y$ given by
$${1\over J(\kappa )}\psi _{-\kappa }(y)\mu _0(dy)$$
converges to the Dirac measure $\delta $ at the base point $e$ of
$\mathcal Y$ in the sense of tight convergence of measures. The
proof of this fact is based on the following lemma.

\begin{lem}\label{suitedelta}
Let $\mathcal K$ be a compact topological space, and $\mu $ a
positive measure on it such that every non empty open set has a
positive measure. Let $q\geq 0$ be a continuous function on
$\mathcal K$ which attains its maximum at only one point $x_0$.
Define, for $n\in \mathbb N$,
$$a_n=\int _{\mathcal K} q(x)^n\mu (dx),$$
and, for a continuous function $\varphi $ on $\mathcal K$,
$$L_n(\varphi )={1\over a_n}\int _{\mathcal K} \varphi (x)q(x)^n\mu (dx).$$
Then
$$\lim _{n\to \infty }L_n(\varphi )=\varphi (x_0).$$
\end{lem}

\emph{Proof} For $0<\alpha <M=\max q$, there exits a constant
$C_{\alpha }$ such that
$$a_n\geq C_{\alpha }\alpha ^n.$$
In fact there is a neighborhood $\mathcal V$ of $x_0$ such that
$q(x)\geq \alpha $ for $x\in \mathcal V$, and
$$a_n\geq \mu (\mathcal V)\alpha ^n.$$
Let $\mathcal W$ be a neighborhood of $x_0$. For $x\in {\mathcal
K}\setminus \mathcal W$, $q(x)\geq \beta <M$. Choose $\alpha $ such
that $\beta <\alpha <M$. then
$${1\over a_n}\int _{\mathcal K\setminus \mathcal W}q(x)^n\mu (dx)
\leq {1\over C_{\alpha }} \mu (\mathcal K)\Bigl({\beta \over \alpha
}\Bigr)^n,$$ and
$$\lim _{n\to \infty }{1\over a_n}\int _{\mathcal K\setminus \mathcal W}q(x)^n\mu
(dx)=0. \qquad\qquad\qquad \Box $$

\bigskip

The function $\psi _{-\kappa }$ attains its maximum $M=1$ only at
$y=e$. Therefore Lemma (\ref{suitedelta}) applies. It follows thet,
for every $\bf m$,
$$\lim _{\kappa \to \infty }{1\over J(\kappa )}a_{\kappa }({\bf m})=1.$$

\medskip

(b) Assume that $\kappa \in \mathbb N$. From the Bernstein identity
(\ref{bernstein}) it follows that the spherical Fourier coefficients
$a_{\kappa}({{\bf m}})$ satisfy the following recursion relation :
$$
\gamma_{-\kappa}(\rho-{{\bf m}})a_{\kappa}({{\bf
m}})=b(-\kappa)a_{\kappa-1}({{\bf m}}).
$$
Since, for ${\bf m}= 0$,
$\gamma_{-\kappa}(\rho)a_{\kappa}(0)=b(-\kappa)a_{\kappa-1}({\bf
0})$ we obtain
$$
\frac{\gamma_{-\kappa}(\rho-{{\bf m}})}{\gamma_{-\kappa}(\rho)}
\frac{a_{\kappa}({{\bf m}})}{a_{\kappa}(0)}
=\frac{a_{\kappa-1}({{\bf m}})}{a_{\kappa-1}(0)}.
$$
Furthermore, since $a_{\kappa }(0)=J(\kappa )$, by (a)
$$
\lim_{\kappa\to\infty}\frac{a_{\kappa}({{\bf m}})}{a_{\kappa}(0)}=1.
$$
The sequence given by
$$
\phi(\kappa)=\frac{\prod_{j=1}^{n}\Gamma(\kappa+\frac12+\delta-\rho_j)\Gamma(\kappa+\frac12+\delta+\rho_j)}
{\prod_{j=1}^{n}\Gamma(\kappa+\frac12+\delta-\rho_j+m_j)\Gamma(\kappa+\frac12+\delta+\rho_j-m_j)},
$$
satisfies the same recursion relation as $\frac{a_{\kappa}({{\bf
m}})}{a_{\kappa}(0)}$ does. Moreover, from the asymptotic
equivalence
$$\frac{\Gamma(z+a)}{\Gamma(z+b)}\sim
z^{a-b}\left(1+\frac1{2z}(b-a)(b+a-1)+O(z^{-2})\right),\quad{\rm
as}\:{z\to\infty},
$$
 it follows that
 $$
 \lim_{\kappa\to\infty}\phi(\kappa)=1.
 $$
 Therefore we have proved the theorem for $\kappa\in \mathbb N$.

\medskip

\noindent(c) By using Carlson's theorem, as we did in the proof of
Theorem (\ref{thmhua}), we conclude that the statement is
 still valid for $\Re\kappa\geq0.\qquad\Box$
\medskip

\subsection{Restriction of the representation $\pi_s$ to the compact
group $U$} The compact symmetric space $\mathcal Y=U/K$ is also of
Makarevich type. Indeed,
$$
U=\{g\in\Conf(V)\:|\:\alpha\circ j\circ g\circ j\circ\alpha=g\}_o.
$$
From the generalized cross-ratio (\ref{cross-ratio}) we define the
kernel $F^c(u,v)$ by
\begin{eqnarray*}
  F^c(u,v) &=& \{u,v,-\alpha(u^{-1}),-\alpha(v^{-1})\} \\
  &=&{\Delta \bigl(u+\alpha (u^{-1})\bigr)\Delta \bigl(v+\alpha
(v^{-1})\bigr) \over \Delta \bigl(u+\alpha (v^{-1})\bigr)\Delta
\bigl(v+\alpha (u^{-1})\bigr)}
\end{eqnarray*}
and $F^c(u,0)=h(u,-u)$.

Similarly to what we did in Section 2.2 we twist the bilinear form
$\mathfrak B_s$ by the involution $j\circ\alpha$, in other words we
replace $y$ by $\alpha(y^{-1})$, and by introducing multiplication
operators
$$
M^c_sf(x)=h(x,-x)^{s+\frac n{2r}}f(x)
$$
we eventually obtain a $U$-invariant bilinear form on
$C^{\infty}(\mathcal Y)$:
$$
\mathfrak{B}_s^c(F_1,F_2)=\iint_{\mathcal{Y}\times\mathcal{Y}}B_s^c(x,y)F_1(x)F_2(y)\mu_0(dx)\mu_0(dy),
$$
where $B_s^c(x,y)=F^c(x,y)^{\frac {r_0}r \nu}$.
\medskip

The spherical dual $\hat U_K$ is parameterized by the set of weights
${\bf m}=(m_1,\ldots,m_n)$ described in the previous section. To
such a weight $\bf m$ corresponds a (class of) unitary spherical
representation $\Pi _{\bf m}$ on a finite dimensional vector space
$\mathcal H_{\bf m}$ of dimension $d_{\bf m}$. The highest weight
$\mu $ of $\Pi _{\bf m}$ is given by
$$
\langle \mu,\xi \rangle =-\sum _{j=1}^{r_0} m_jt_j,$$ if $\xi
(x)=\sum t_jL(c_j)x$. In the space $\mathcal H_{\bf m}$ there is a
normalized $K$-fixed vector $v_{\bf m}$. The Fourier coefficient
$\hat f(\bf m)$ of an integrable function $f$ on $\mathcal{Y}$ is
the vector in $\mathcal H_{\bf m}$ defined by
$$
\hat f({\bf m})=\int_{U/K}\Pi _{\bf m}(g)v_{\bf m}f(g)dg.
$$
The map $f\mapsto\hat f(\bf m)$ intertwines the left regular
representation of $U$ and $\Pi _{\bf m}$.

\begin{thm}
For  a $\kappa\in\mathbb N$ the Berezin form $\mathfrak B_s^c$ is
positive definite. Therefore the restriction of $\pi_s$ to the
compact group $U$ decomposes into a direct sum of spherical
principal series representations $\Pi_{\bf m}$ of $U$ according to
the following Plancherel formula
$$
\mathfrak B_s^c(f,\bar f)=\sum_{{\bf m}\in\hat U_K}d_{\bf
m}a_{\kappa}({\bf m}) \Vert\hat f({\bf m})\Vert^2,
$$
where $\kappa=\frac r{r_0}(s-\frac{n}{2r})$, $f\in C(U/K)$.
\end{thm}

By the first part of the proof of Theorem 2.8, in some sense, the
canonical representations of the compact dual of $G$ tend to the
left regular representation when the parameter $\kappa$ goes through
the negative integer points. In fact, for $f\in \mathcal C_c
(\mathcal Y)$,
$$
\lim_{\kappa\to\infty} \frac1{J(\kappa )}\int _{\mathcal Y\times
\mathcal Y} \mathfrak B_{\kappa }^c(x,y)f(x)f(y)\mu _0(dx)\mu _0(dy)
=\int _{\mathcal Y}|f(y)|^2\mu _0(dy).
$$
 Similarly, in the non compact case, as
$\nu\to \infty $, the probability measure $ I(\nu)^{-1}\psi _{\nu
}(x)\mu(dx)$ converges to the Dirac measure $\delta _e$ at the
identity element $e$. It follows that, for $f\in \mathcal C
_c(\mathcal X)$,
$$\lim _{\nu\to \infty }
{1\over I({\nu })}\int _{{\mathcal X}\times {\mathcal X}} {\mathcal
B}^{\alpha}_{\nu }(x,y)f(x)\overline{f(y)}\mu(dx)\mu(dy) =\int
_{\mathcal X} |f(x)|^2\mu(dx).$$ In some sense, as $\nu \to \infty
$, the canonical representation $T_\nu$ tends to the regular
representation of $G$ on $L^2({\mathcal X})$.

In conclusion we should mention that originally the one parameter
family of Berezin kernels was introduced in the context of
quantization of K\"ahlerian symmetric spaces in \cite{[Ber]}. The
last remark says that the Berezin symbolic calculus "deforms" the
usual point-wise multiplication of functions. We refer the reader to
\cite{[AO],[UU]} for more details on asymptotics of the Berezin
transform.
\bibliographystyle{amsplain}

\end{document}